\documentclass[12pt]{article}
\scrollmode

\usepackage{color}

\usepackage{amsfonts}
\usepackage{amsmath}
\usepackage{amsthm}
\usepackage{amssymb}
\usepackage{graphicx}
\usepackage{enumitem}

\allowdisplaybreaks

\title{A stronger topology for the Brownian web\footnote{In celebration of Chuck Newman's $70^{\mbox{\scriptsize th}}$ birthday.}} 
\date{}
\author{L.~R.~G.~Fontes\footnote{IME-USP, Rua do Mat\~ao 1010, 05508-090, S\~ao Paulo SP,  Brazil, lrfontes@usp.br}
\thanks{Partially supported by CNPq grant 311257/2014-3 and FAPESP grant 2017/10555-0.}
}

\newtheorem{theo}{Theorem}

\newtheorem{defin}[theo]{Definition}
\newtheorem{lem}[theo]{Lemma}
\newtheorem{rmk}[theo]{Remark}

\def\Z{{\mathbb Z}}
\def\R{{\mathbb R}}

\def\cc{{\cal C}}
\def\cd{{\cal D}}

\def\ch{{\cal H}}

\def\ci{{\cal I}}

\def\cm{{\cal M}}
\def\cn{{\cal N}}

\def\cp{{\cal P}}

\def\cR{{\cal R}}
\def\cs{{\cal S}}

\def\T{{\cal T}}

\def\cw{{\cal W}}
\def\bcw{\bar{\cal W}}

\def\d{\delta}
\def\D{\Delta}

\def\eps{\epsilon}
\def\teps{\tilde\epsilon}

\def\lf{\lfloor}
\def\rf{\rfloor}

\def\s{\sigma}
\def\a{\alpha}

\def\G{\Gamma}
\def\l{\lambda}

\def\eps{\varepsilon}

\def\tt{\tilde t}
\def\tp{\tilde p}
\def\tq{\tilde q}

\def\bd{\bar d}
\def\hd{\tilde d}
\def\dc{\bar d_\ch}
\def\hh{\ch_{\cc}}
\def\hhz{\ch^0_{\cc}}
\def\hhl{\ch^L_{\cc}}
\def\hhhz{\hat\ch^0_{\cc}}

\def\hdh{d_{\hh}}

\def\dch{d_{\ch}}

\def\mud{\mu^{(\d)}}
\def\bud{\bar\mu^{(\d)}}

\def\yd{Y^{(\d)}}
\def\byd{\bar Y^{(\d)}}
\def\pd{\cp^{(\d)}}
\def\bpd{\bar\cp^{(\d)}}
\def\pdz{\cp^{(\d)}_0}
\def\pdl{\cp^{(\d)}_L}
\def\cpz{{\cal P}_0}
\def\cpl{{\cal P}_L}
\def\pdze{\cp^{(\d)}_{0,\eps}}
\def\pdzz{\cp^{(0)}_{0,\eps}}
\def\pze{\cp_{0,\eps}}

\def\crd{\cR^{(\d)}}
\def\ctd{\T^{(\d)}}
\def\td{T^{(\d)}}
\def\tt{\bar t}
\def\tf{\bar f}
\def\hf{\tilde f}

\def\={&=&}

\def\nn{\nonumber}

\begin{document}

\maketitle

	\begin{abstract}
		
	We propose a metric space of {\em coalescing} pairs of paths on which we are able 
	to prove (more or less) directly convergence of objects such as the {\em persistence probability} 
	in the (one dimensional,nearest neighbor, symmetric) voter model
        or the diffusively rescaled {\em weight distribution} in a silo model 
        (as well as the equivalent {\em output distribution} in a river basin model), 
        interpreted in terms of (dual) diffusively rescaled coalescing random walks, 
	to corresponding objects defined in terms of the Brownian web.

	\end{abstract}
	
	
	\smallskip
	


\section{Introduction}
\label{intro}

\setcounter{equation}{0}

Chuck Newman first noticed that the convergence results of~\cite{kn:finr} were not by themselves enough to show convergence of
the {\em persistence probability} in the voter model to a corresponding probability in the Brownian web. 
As is well known, the voter model may be described in terms of a system of coalescing random walks, which in turn
has in~\cite{kn:finr} been shown to converge weakly to the Brownian web (under diffusive scaling). 
As pointed out by Chuck, the issue has to do with the topology in the space of trajectories adopted in~\cite{kn:finr}.
Let us elaborate these points, and introduce other instances where they come up.

\paragraph{Persistence in the voter model.}

Let us suppose that individuals are placed on the sites of $\Z$, and that initially all the individuals on $2\Z$ have each a 
different opinion on a certain matter. At 
time $n\geq1$, each individual of $2\Z+(n\!\!\mod 2)$\footnote{That is, the even or odd sublattices of $\Z$, depending on the parity of $n$, respectively.} adopts the
opinion  at time $n-1$ of a nearest neighbor individual (on $\Z$) chosen uniformly at random, independently from other individuals and choices.

At time $n=2k$, the opinion of the individual at the origin may be traced back to the opinion of an individual at time $0$. Let the choice of neighbor by individual $i$ at time $\ell$, as described above, be graphically represented as an edge from $(i,\ell)$ to either $(i-1,\ell-1)$ or $(i+1,\ell-1)$, depending whether the chosen neighbor was the one to the left or the one on the right, respectively. Let us orient each such edge {\em backwardly} from $\ell$ to $\ell-1$. We thus obtain a system of {\em backward} paths, one from each point of $\Z^2_{e,+}$, the even sublattice of $\Z^2$ in the upper half plane down to a point of $2\Z\times\{0\}$, given by a simple random walk path from the initial point to the final point.
We say that the path model thus obtained is {\em dual} to the voter model.
The paths of the path model are independent before meeting, and when they meet, they coalesce into the single path starting at the meeting space-time point (with time running backwards).

The collection of diffusively rescaled such coalescing backward paths converges as the scaling parameter vanishes to the Brownian web. This is the convergence result of~\cite{kn:finr} mentioned above.

Let us now define the persistence probability of the voter model as the probability that the individual at the origin does not change opinion from time $m=2j$ to time $n$. This object was considered in~\cite{kn:fins} for the equivalent zero temperature Metropolis dynamics for the one dimensional nearest neighbor Ising model. In terms of the dual path model, it is equivalent to the probability that all backward paths starting from $\{(0,2\ell);\,\ell=j,\ldots,k\}$ coalesce before or at time $0$.

Rescaling the model diffusively, with $m$ and $n$ also rescaled accordingly, and in such a way that $m/n\sim\a\in(0,1)$, 
it is natural to expect that the persistence probability for the voter model will converge as the scaling parameter vanishes to
the probability that all the (forward) paths from the Brownian web starting from the vertical segment $\{0\}\times[0,\a]$ coalesce before time $1$. But this does not follow from the convergence result of~\cite{kn:finr} by itself.

Before proceeding with this particular point, let us discuss another model, where a similar issue comes up.

\paragraph{Weight distribution in a silo model.}  

In~\cite{kn:clmnw} a silo model is proposed as follows. A unit weight bead is located at every point of $\Z^2_{e,+}$. Each bead, located at say $(i,\ell)$ supports the full weight of either the one at $(i-1,\ell+1)$ or the one at $(i+1,\ell+1)$, which one is determined by a uniformly random choice, independent from bead to bead.
The {\em total weight} supported by a bead is the sum of the weights of all the beads it supports from all the rows of beads above it, including its own weight.

In order to understand this quantity,
let us start by putting downward directed edges from the supported to the supporting beads, and consider the downward paths starting from each bead to a bead at $2\Z\times\{0\}$ it connects to via those downward edges. Each pair of these paths are independent random walk paths till they meet, when they coalesce into a single random walk path.

Let us now consider the following set of dual upward directed edges in $\Z^2_{o,+}$, the {\em odd} sublattice of $\Z^2$ in the upper half plane. An upward edge from $(j,k)$ in $\Z^2_{o,+}$ to either $(j-1,k+1)$ or $(j+1,k+1)$ is present if and only if it does not cross any of the downward edges described in the previous paragraph. Let us now consider the upward paths starting from the points of $\Z^2_{o,+}$ and going up along the upward edges of $\Z^2_{o,+}$ just described.
One may readily check that any two distinct such upward paths are independent until they meet, after which time they coalesce.
We say that the family of upward paths is {\em dual} to the system of downward paths\footnote{This is a different notion of duality than the one we described for the voter model above.}.

It may be also checked that the total weight supported by a bead located at $(i,\ell)$ is given by the number of beads enclosed by the upward paths starting from $(i-1,\ell)$ and $(i+1,\ell)$.

We want to study the total weight distribution at the bottom of the silo, which may be formulated as the measure on 
$2\Z$ assigning to $\{i\}$ the total weight supported by the bead at $(i,0)$.

By properly rescaling this measure, we may expect that it converges weakly to the measure on $\R$ that assigns to a 
finite interval $[a,b]$ the area encompassed by the following paths from the Brownian web: the leftmost path from $(a,0)$
and the rightmost path from $(b,0)$\footnote{In the Brownian web we may find multiple paths starting from a single space-time point.}. Again this does not follow from the convergence result of~\cite{kn:finr}.

\paragraph{Output distribution in a river basin model.} 

The river basin/drainage system proposed in~\cite{kn:s} 
(see also~\cite{kn:rr} for a broader discussion of such models) may be described as follows. From each site $(i,\ell)$ of $\Z^2_{e,+}$ there flows, along the appropriate edge, a unit volume of water to either $(i-1,\ell+1)$ or $(i+1,\ell+1)$, which one is determined by a uniformly random choice, independently from one originating site to the other. The downward flows to $2\Z\times\{0\}$ are thus described as downward random walk paths from the sites of $\Z^2_{e,+}$ to sites of $2\Z\times\{0\}$,
abstractly identical to the ones from the silo model discussed above. We are interested in the distribution of total water flow on the sites of $2\Z\times\{0\}$. This is abstractly identical to the weight distribution measure introduced in the previous example, and thus should converge also as before to the Brownian web measure introduced above. 

Let us next discuss why the convergence result of~\cite{kn:finr} is not enough to establish the convergences in the above examples.

\paragraph{The trouble with the topology of~\cite{kn:finr}.} It comes up already when considering two coalescing paths, say starting from different space-time points (each in the segment the vertical segment $\{0\}\times[0,\a]$, as in the first example; or both on $\R\times\{0\}$, as in the other examples); these paths should not touch till a {\em coalescence time}, up from which they are identical. Precisely the coalescence time should be continuous in the sup kind of metric of~\cite{kn:finr} in order for
the convergences in the examples to follow from the convergence result of~\cite{kn:finr}. But it is not: two pairs of coalescing paths may be close in sup metric without their coalescence times being close. Indeed, it is not continuous at any pair of coalescing paths; even the {\em coalescence property} (as implicitly given above) is not continuous in the metric space where 
the Brownian web is defined in~\cite{kn:finr}.

\paragraph{A new metric space.} 

This paper is an effort to address the issue raised above by introducing a new metric space for the Brownian web which takes into consideration the coalescence property and its continuity (in an appropriate sense). We discuss the new space in Section~\ref{mod}, and in Section~\ref{bw} show that an appropriate description of the Brownian web belongs to it. There is {\em not} an effort at generality, so this space will accommodate the simple coalescing random walks of the examples, but not non simple ones. 
It is possible that a mild variation of this space and metric will include the latter models, but we do not discuss the matter further here.

We then show in Section~\ref{rw} a weak convergence result of rescaled random walk paths starting from all space locations {\em at fixed times} to the corresponding Brownian web paths in this new space. This is enough to establish the convergences of the examples more or less directly in Section~\ref{ex}. Unfortunately, we could not show weak convergence of the full family of paths starting from {\em all} space-time points to the full Brownian web in this new  space; more on this point in the last section of the paper.


\section{Topological set-up}\label{mod}

\setcounter{equation}{0}

For simplicity, we will consider paths starting from a bounded region of the space-time plane 
--- which we take to be the rectangle $\cR:=[-1,1]\times[0,1]$ ---, rather than the compactified 
$\R^2$ as in~\cite{kn:finr}. We will actually consider pairs of paths starting from $\cR$: we 
start with a space of pairs of coalescing paths; and then a (Hausdorff) space of compact sets 
of the former space. The latter space will be the sample space of the (restricted) Brownian web.


\paragraph{Pairs of coalescing paths.}

The first space contains pairs of paths from a single path space which we define now. Let
\begin{equation}
\Pi'=\Pi'_{\cR}=\bigcup_{t_0\in[0,1]}C[t_0]\times\{t_0\},
\end{equation}
where $C[t_0]$ denotes the set of continuous functions $f:[t_0,\infty)\to\R$ such that 
$f(t_0)\in[-1,1]$. 
The elements of $\Pi'$ should be seen as continuous paths in the
space-time plane starting in $\cR$, with time running upwards, each such path starting at a given time $t_0$ from $f(t_0)$; 
such an element is denoted $(f,t_0)$. 

To describe the first space, we need a notion of coalescence of a pair of paths from $\Pi'$, to be defined now.
\begin{defin}\label{coal}
For $i=1,2$ let $t_i\in[-1,1]$, and suppose $(f_i,t_i)\in\Pi'$. Let 
\begin{equation}\label{tc}
t^c=\inf\{t>t^+:=t_1\vee t_2:\,f_1(t)=f_2(t)\},
\end{equation}
with the usual convention that $\inf\emptyset=\infty$.
We will say that $(f_1,t_1)$ and $(f_2,t_2)$ {\em coalesce}, or are {\em coalescing}, 
or is a {\em coalescing pair (of paths)} if either $t^c=\infty$ or 
$f_1(t)=f_2(t)$ for $t\geq t^c$.
\end{defin}

\begin{rmk}\label{coar}\mbox{ }
	
\begin{enumerate}
    \item The {\em coalescence time} $t^c$ may equal $t^+$ --- a case that happens only if $f_1(t^+)=f_2(t^+)$.
          In particular, two identical paths of $\Pi'$ are coalescing. But it may happen that $f_1(t^+)=f_2(t^+)$ and $t^c>t^+$; this indeed takes place in the Brownian web.
    \item In each coalescing pair of paths, we may distinguish a {\em left path} and a {\em right path}, 
          such that the left path is always to the left of (or coincides with) the right path (for $t\geq t^+$)
          --- that is, $f_1(t)\leq f_2(t)$ for $t\geq t^+$. 
          We will denote a coalescing pair by $[(f_1,t_1),(f_2,t_2)]$, where the first path, namely $(f_1,t_1)$, will always
          denote the left path; occasionally we will also use the notation $[f_1,f_2]$ for short. In this context, we
          will occasionally write about an {\em ordered} pair of paths
\end{enumerate}	
\end{rmk}	

The second space is finally defined as
\begin{equation}\label{coas}
\cc=\{(f_1,t_1),(f_2,t_2)\in\Pi':\,(f_1,t_1)\mbox{ and }(f_2,t_2)\mbox{ coalesce}\}.
\end{equation}

We now want to put a metric in $\cc$ under which it is complete and separable. Completeness is the tricky
thing. Simple variants of the sup metric do not prevent Cauchy sequences which, so to say, want to converge
to pairs of paths of $\Pi'$ which touch and later coalesce, and are thus not in $\cc$. In order to have such sequences not 
be Cauchy, we add to (a suitable variant of) the sup metric a term tailored for this end. We discuss that term next. 

\paragraph{Sup metrics in $\Pi'$ and $\cc$.}

Let us first consider the sup metric in $\Pi'$ defined as follows: given $(f,t_0)$, $(g,s_0)\in\Pi'$, let
\begin{equation}
\label{d}
d'((f,t_0),(g,s_0))=
|t_0-s_0|\vee\sum_{n\geq1}2^{-n}\left(\sup_{t\in[0,n]}|\hat{f}(t)-\hat{g}(t)|\wedge1\right),
\end{equation}
where, given $(f,t_0)\in\Pi'$, $\hat f:[0,\infty]\to\R$ 
is such that $\hat f(t)=f(t)$ for $t\geq t_0$, and $\hat f(t)=f(t_0)$ for $t\in[0, t_0]$.

Let $\bar d$ denote the Hausdorff-type metric induced in $\cc$ by the metric $d'$ in $\Pi'$, namely, 
given two pairs of coalescing paths $[(f_1,t_1),(f_2,t_2)]$ and $[(g_1,s_1),(g_2,s_2)]$
(each pair in $\cc$), let
\begin{equation}
\label{bd}
\bar d([(f_1,t_1),(f_2,t_2)],[(g_1,s_1),(g_2,s_2)])=
\max_{i=1,2}\min_{j=1,2}d'((f_i,t_i),(g_j,s_j)).
\end{equation}

It may be readily checked that $\bar d$ is a metric in $\cc$, but {\em not} a complete one: there
are Cauchy sequences $[(f^{(n)}_1,t^{(n)}_1),(f^{(n)}_2,t^{(n)}_2)]$ in $(\cc,\bar d)$ such that as $n\to\infty$
$(f^{(n)}_1,t^{(n)}_1)\to(g_1,s_1)$ and $(f^{(n)}_2,t^{(n)}_2)\to(g_2,s_2)$ in $(\Pi',d')$, with 
$(g_1,s_1)$ and $(g_2,s_2)$ such that there exist
$t',s^*$ with $t'\in(s^+,s^*)$ satisfying that $g_1(t)\ne g_2(t)$ for $t\in(s^+,t')\cup(t',s^*)$, 
and $g_1(t)=g_2(t)$ for $t=t',s^*$; thus $[(g_1,s_1),(g_2,s_2)]\notin\cc$. (Another way to put it is that 
$\cc$ is not closed in $(\Pi'\times\Pi',\bar d)$.)

We will then add a term to $\bar d$ in order to prevent this situation and have a resulting complete (and separable)
metric in $\cc$, as follows.

\paragraph{Pods.} We rely on the following concept of a {\em pod} of a pair of coalescing paths in $\cc$, defined roughly 
as the portion of the pair of paths between the beginning of the later path and the coalescence point. 
We find it necessary to separate the portions of the paths close to the beginning of the pod (the beginning of the later path) from the portions close to the end of the pod (the coalescence time)
thus producing a {\em standard pod}. 

For convenience, things will be defined in terms of the difference between the paths, as follows.
%
Given a pair $[(f_1,t_1),(f_2,t_2)]\in\cc$, recall that $t^+=t_1\vee t_2$ and $t_c$ is the coalescence point. 
Let  $\Psi(\cdot)=\tanh(\cdot)$ and make
\begin{equation}\label{tt}
\tilde t^+= \Psi(t^+),\,\tilde t^c=\Psi(t^c),
\end{equation}

%

If $t^+<t^c$, then let $\tp:[\tilde t^+,\tilde t^c)\to(0,\infty)$ be such that 
\begin{eqnarray}\label{pp}
\tp(x)=f_2(\Psi^{-1}(x))-f_1(\Psi^{-1}(x))
\end{eqnarray}
$x\in[\tilde t^+,\tilde t^c)$, where $\Psi^{-1}$ is the inverse function of $\Psi$. 
We recall that, according to the convention set above, $f_2\geq f_1$.

If $t^+=t^c$, then let $\tp:\{\tilde t^+\}\to\{0\}$. 

We note that if $t^c<\infty$, then $\tp(t^c):=\lim_{x\uparrow t^c}\tp(x)=0$ in all cases.

\begin{defin}\label{def:pod}
We call $\tp$ the {\em pod (function)} of $[(f_1,t_1),(f_2,t_2)]$. We will refer below to $t^+$ and $t^c$ as
the beginning and end of the pod, respectively. We will also refer to the {\em length} of the pod, meaning 
$t^c-t^+$, and also to the {\em diameter} of the pod, meaning $\sup_{x\in(t^+,t^c)}\tp(x)$.
We finally define the {\em dimension} of a pod as the max of its length and diameter.
\end{defin}
Now let $\tilde t^m$ be the midpoint of $[\tilde t^+,\tilde t^c]$ and for $x\in[0,1)$ make
\begin{equation}\label{p}
p(x)=
\begin{cases}
\tp(x+\tilde t^+),\,& x\in[0,\tilde t^m-\tilde t^+];\\
\tp(x-(1-\tilde t^c)),\,& x\in[1-(\tilde t^c-\tilde t^m),1);\\
\tp(\tilde t^m)+\frac12[1-(\tilde t^c-{\tilde t}^m)],\,& x=1/2.
\end{cases}
\end{equation}
If $t^+=t^c$, we define $p(1)=0$. We then complete the definition of $p$ in $[0,1)$ by
making it linear in $[\tilde t^m-\tt^+,1/2]$ and in $[1/2,1-(\tilde t^c-\tilde t^m)]$. 

We call $p$ the {\em standard pod (function)} of $[(f_1,t_1),(f_2,t_2)]$,
and note that $p$ satisfies
\begin{equation}\label{prop}
p:(0,1)\to(0,\infty); \,\lim_{x\uparrow1}p(x)=0,\mbox{ if }t^c<\infty;\,p(0)\geq0,
\end{equation}
with equality if and only if $f_1(t^+)=f_2(t^+)$. Additionally, $p$ has derivative $1$ and $-1$
in $(\tilde t^m-\tilde t^+,1/2)$ and $(1/2,1-(\tilde t^c-\tilde t^m))$, respectively.


For $p,q$ with the above properties (namely,~(\ref{prop}) and continuity) let
\begin{equation}\label{Delta}
\Delta(p,q)=\sum_{n\geq2}2^{-n}\left(\sup_{x\in\left[\frac1n,1-\frac1n\right]}\left| \frac1{p(x)}-\frac1{q(x)}\right| \wedge1\right).
\end{equation}

The following result will be needed in Section~\ref{rw} below.

\begin{lem}\label{lem:pod}
	If $\tp$ and $\tq$ are pods with dimension at most $\s$, then $\Delta(p,q)\leq\varphi(\s)$,
	with $\lim_{\s\to0}\varphi(\s)=0$, where $p$ and $q$ are the respective standard pods of
	$\tp$ and $\tq$.
\end{lem}	

\noindent{\bf Proof.} It follows from~(\ref{p}) that, for all $\s\leq1/2$, $p$ and $q$ are both linear 
separately in the intervals $[\s,1/2]$ and $[1/2,1-\s]$, with inclinations $1$ and $-1$ in each interval, 
respectively. This implies that 
\begin{equation}\label{bdd}
\sup_{x\in[\s^{1/3},1-\s^{1/3}]}\left| \frac1{p(x)}-\frac1{q(x)}\right|\leq\frac{2\s}{\s^{2/3}-\s^{2}}
=:\varphi'(\s).
\end{equation}
It readily follows that $\varphi(\s)=\varphi'(\s)+\sum_{n\geq\s^{-1/3}}2^{-n}$ verifies the claim.
$\square$

\paragraph{A stronger metric in $\cc$.}

We discuss now the term to be added to $\bar d$ in order to have a metric in $\cc$ making it complete. 

Let  $[(f_1,t_1),(f_2,t_2)]$ and $[(g_1,s_1),(g_2,s_2)]$ be two coalescing pairs of $\cc$, 
and let $p,q$ be the standard pods respectively of $[(f_1,t_1),(f_2,t_2)],[(g_1,s_1),(g_2,s_2)]$. 


Let us define
\begin{eqnarray}\nn
&\hd([(f_1,t_1),(f_2,t_2)],[(g_1,s_1),(g_2,s_2)])&\\\label{met}
&= \bd([(f_1,t_1),(f_2,t_2)],[(g_1,s_1),(g_2,s_2)]) + \D(p,q)&
\end{eqnarray}
for $[(f_1,t_1),(f_2,t_2)],\,[(g_1,s_1),(g_2,s_2)]\in\cc$.


\paragraph{$(\cc,\hd)$ is a complete and separable metric space.}

That $\hd$ is a metric in $\cc$ is quite clear; separability is also fairly clear after a moment's thought.

To argue completeness (sketchily), let $[(f^{(n)}_1,t^{(n)}_1),(f^{(n)}_2,t^{(n)}_2)]$, $n\geq1$, be a
Cauchy sequence in $(\cc,\hd)$. 
Let $t^+_n$ and $t^c_n$ be the beginning and end of the respective pods.

From the $\bd$ part of $\hd$, we have that, for $i=1,2$, there exist 
$(f_i,t_i)\in\Pi'$ such that, as $n\to\infty$, 
\begin{eqnarray}\label{conv1}
&(f^{(n)}_i,t^{(n)}_i)\to(f_i,t_i)\mbox{ in }(\Pi',d');&\\\label{conv2}
&t^{(n)}_i\to t_i,\, t^+_n\to t^+=t_1\vee t_2\mbox{ in }\R.&
\end{eqnarray}
And combining this with the $\D$ part of $\hd$, we find that 
\begin{equation}\label{conv1a}
t^c_n\to t^c\mbox{ in }\bar\R
\end{equation}
for some $t^c\in\bar\R:=\R\cup\{\infty\}$. 
It may be readily checked that if $t^c<\infty$, then $f_1(t)=f_2(t)$ for $t\geq t^c$.

If $t^+=t^c$, then we immediately have that $(f_1,t_1)=(f_2,t_2)$, and 
then of course $[(f_1,t_1),(f_2,t_2)]\in\cc$ in this case.

Let us now verify that if $t^+<t^c$, then 
\begin{equation}\label{podp}
f_1(t)\ne f_2(t)\mbox{ for }t\in(t^+, t^c),
\end{equation}
and thus we have that 
$[(f_1,t_1),(f_2,t_2)]\in\cc$ also in this case.

From~(\ref{conv2}) and~(\ref{conv1a}), we have that  $\tilde t^+_n\to\tilde t^+$ and $\tilde t^c_n\to\tilde t^c$, and $t^+<t^c$ implies that
$\tilde t^+<\tilde t^c$. Given now $[s,s']\subset(t^+, t^c)$, we have that $\tilde s=\Psi(s),\tilde s'=\Psi(s')$ are such
that $[\tilde s,\tilde s']\subset(\tilde t^+, \tilde t^c)$. 

We will at this point assume for simplicity, but without compromising the generality of conclusion of the argument, 
that $\tilde t^m\in(\tilde s,\tilde s')$.

Let $p_n$, $n\geq1$, be the sequence of standard pods of $[(f^{(n)}_1,t^{(n)}_1),(f^{(n)}_2,t^{(n)}_2)]$. Since by
hypothesis $\D(p_m,p_n)\to0$ as $m,n\to\infty$, it follows that there exists $p:(0,1)\to(0,\infty)$ such that
$\frac1{p_n}\to\frac1{p}$ uniformly over closed intervals of $(0,1)$. 
Notice also that from the $\bd$ part of the metric, $p$ must be of the form~(\ref{p},\ref{pp}).
It follows that $p>0$ in $[\tilde s-\tt^+,1-\tilde s']$; in turn $\tp>0$ in $[\tilde s,\tilde s']$,
from which it immediately follows that $f_1\ne f_2$ in $[s,s']$, and~(\ref{podp}) is established. 

To finish the completeness argument, it remains to show that 
\begin{equation}\label{comp}
[(f^{(n)}_1,t^{(n)}_1),(f^{(n)}_2,t^{(n)}_2)]\to[(f_1,t_1),(f_2,t_2)]
\end{equation}
in $(\cc,\hd)$ as $n\to\infty$, but this follows readily from~(\ref{conv1}-\ref{conv1a}) and the arguments of
the previous paragraph.

\begin{rmk}\label{rmk:conv}
	If $[(f_1,t_1),(f_2,t_2)]$ and $[(f^{(n)}_1,t^{(n)}_1),(f^{(n)}_2,t^{(n)}_2)]$, $n\geq1$,
	are in $\cc$ and~(\ref{conv1},\ref{conv1a}) hold with $t^c$ the coalescence time of $[(f_1,t_1),(f_2,t_2)]$, 
	then, as may be readily checked,
	\begin{equation}\label{convh}
	[(f^{(n)}_1,t^{(n)}_1),(f^{(n)}_2,t^{(n)}_2)]\to[(f_1,t_1),(f_2,t_2)]\mbox{ in }(\cc,\hd).
	\end{equation} 
\end{rmk}

\paragraph{Hausdorff metric space.} We may now define $(\hh,\hdh)$ as the Hausdorff metric
space associated to $(\cc,\hd)$, namely
\begin{equation}\label{hh}
\hh=\{K\subset\cc:\,K\mbox{ is non empty and compact}\}
\end{equation}
and, given $K,K'\in\hh$,
\begin{eqnarray}\nn
&\hdh(K,K')=&\\\label{hdh}
&{\displaystyle\sup_{[f_1,f_2]\in K}\inf_{[g_1,g_2]\in K'}\hd([f_1,f_2],[g_1,g_2])\vee
\sup_{[g_1,g_2]\in K'}\inf_{[f_1,f_2]\in K}\hd([f_1,f_2],[g_1,g_2])}&
\end{eqnarray}
(see the last sentence of the 4th point of Remark~\ref{coar} above).

\begin{rmk}
Given the properties of completeness and separability of $(\cc,\hd)$, established above, 
it follows that $(\hh,\hdh)$ has the same properties.
\end{rmk}


\section{The Brownian web as a compact set of coalescing pairs of paths}
\label{bw}

\setcounter{equation}{0}

The Brownian web was defined in~\cite{kn:finr} as a closed set of coalescing Brownian paths 
from $(\Pi,d)$, a single path space similar to $(\Pi',d')$, with the difference that paths 
live in a compactified $\R^2$ --- see Section 3 of~\cite{kn:finr}.
It was then showed to be compact. In this section, we will define a {\em restricted} Brownian web 
as a set of pairs of paths from $(\cc,\hd)$, and then show that it is compact,
and thus belongs to $\hh$. This will set up the stage for proving 
a (somewhat restricted)
convergence of rescaled coalescing simple random walks to the Brownian web in $(\hh,\hdh)$ in the
next section.

Let us briefly make a definition based on a definition from~\cite{kn:finr}. 
Let $\cd$ be a dense 
countable subset of $\cR$ and let $\cw=\cw(\cd)$ be collection of coalescing 
Brownian paths (with unit diffusion coefficient), one starting from each point
of $\cd$. See~\cite{kn:finr}, Section 3, for details of a construction. Let
now $\bcw=\bcw(\cd)$ be the closure of $\cw$ as a subset of $(\Pi',d')$. 
$\bcw$ is what we call the {\em restricted Brownian web}, and $\cw$ its {\em skeleton}.

It is quite clear that almost surely every pair of paths in $\cw$ is coalescing
according to Definition~\ref{coal} above. $\bcw$ has the same property, as
follows readily from Proposition 4.2 of~\cite{kn:finr}. So the collection of
coalescing (ordered) pairs 
\begin{equation}\label{bpw}
\cp:=\{((f_1,t_1),(f_2,t_2))\in\bcw\times\bcw:\,f_2(t)\geq f_1(t)\mbox{ for }t\geq t_1\vee t_2\}
\end{equation}
is a almost surely subset of $\cc$.

\paragraph{$\cp$ is a compact set of $(\cc,\hd)$.}

We will argue now the almost sure compactness of $\cp$ in $(\cc,\hd)$.

We know from Proposition 3.2 in~\cite{kn:finr} that $\bcw$ is compact in $(\Pi,d)$, so given a 
sequence $[(f^{(n)}_1,t^{(n)}_1),(f^{(n)}_2,t^{(n)}_2)]$, $n\geq1$, from $\cp$, we have that
there exists a subsequence $(n')$ of $(n)$ and $[(f_1,t_1),(f_2,t_2)]\in\cp$ such that
\begin{eqnarray}\label{convbw1}
(f^{(n')}_i,t^{(n')}_i)\to(f_i,t_i)\mbox{ in }(\Pi,d).
\end{eqnarray}
Since Brownian paths do not go to $\pm\infty$ in finite time, we readily conclude that
\begin{eqnarray}\label{convbw2}
(f^{(n')}_i,t^{(n')}_i)\to(f_i,t_i)\mbox{ in }(\Pi',d'),
\end{eqnarray}
and, as pointed out in Remark~\ref{rmk:conv}, 
we will have that
\begin{equation}\label{comp1}
[(f^{(n')}_1,t^{(n')}_1),(f^{(n')}_2,t^{(n')}_2)]\to[(f_1,t_1),(f_2,t_2)]\mbox{ in }(\cc,\hd)
\end{equation} 
once~(\ref{conv1a}) holds. Let us argue the latter point. It follows from the proof of Proposition 4.3 of~\cite{kn:finra}
and~(\ref{convbw1}) that given $\eps>0$, there exist a pair of paths $(f'_1,t'_1),(f'_2,t'_2)\in\cw$ and $n_0$ such that 
for every $n'\geq n_0$, we have that
\begin{eqnarray}\label{id1}
f^{(n')}_1(t)&=&f'_1(t),\,\mbox{ for }t\geq t_1+\eps\\\label{id2}
f^{(n')}_2(t)&=&f'_2(t),\,\mbox{ for }t\geq t_2+\eps,
\end{eqnarray}
and~(\ref{conv1a}) follows. 

\begin{defin}
We will call $\cp$ the {\em Brownian pair web (BPW)}. 
\end{defin}

\begin{rmk}
 As argued above, the BPW is almost surely compact, and as such, it belongs to $(\hh,\hdh)$.
\end{rmk}

\begin{defin}\label{def:ptau}
	For $\tau\in[0,1]$, let 
\begin{equation}\label{ptau}
\cp_\tau=\{[(f_1,t_1),(f_2,t_2)]\in\cp:\,t_1=t_2=\tau\}
\end{equation} 	
be the pairs in $\cp$ both paths of which start from time $\tau$.
\end{defin}

It may be readily checked that $\cp_\tau$ is compact for every $\tau\in[0,1]$ almost surely.


\section{Coalescing simple random walks}
\label{rw}

\setcounter{equation}{0}

Let $Y$ and $\yd$ be respectively the collection of coalescing simple random walk paths, and the
diffusively rescaled collection of the same set of paths, with scale parameter $\d>0$, both collections 
consisting of paths starting from points in the space-time plane. See Introduction 
and Section 6 of~\cite{kn:finr}. We recall that the starting space-time points of $Y$ consist 
of all of $\Z^2$. We will refer below to the starting space-time points of $\yd$ as 
{\em (rescaled) discrete space-time points}.
We will consider the subcollection of paths from $\yd$ starting
from $\cR$. Let us denote that collection by $\byd$, and their starting points by $\crd$. 
Let us now consider the set of (ordered) pairs from $\byd$:
\begin{equation}\label{dpw}
\pd:=\{((f_1,t_1),(f_2,t_2))\in\byd\times\byd:\,f_2(t)\geq f_1(t)\mbox{ for }t\geq t_1\vee t_2\}.
\end{equation}

We call $\pd$ the {\em Discrete Pair Web} (DPW).
We note that, for every fixed $\d>0$, $\pd$ is a finite subset of $\cc$, so it belongs to $(\hh,\hdh)$.


\begin{defin}\label{def:ptaud}
	For $\tau\in[0,1]$, let 
	\begin{equation}\label{ptaud}
	\pd_\tau=\{[(f_1,t_1),(f_2,t_2)]\in\pd:\,t_1=t_2=\d^2\lf\d^{-2}\tau\rf\}
	\end{equation} 	
	be the pairs in $\pd$ both paths of which start from time $\d^2\lf\d^{-2}\tau\rf$.
\end{defin}

We will now argue that as $\d\to0$, the distribution of $\pd_\tau$ in $(\hh,\hdh)$ converges 
$\cp_\tau$ 
defined in~(\ref{ptau}) above.

\paragraph{$\pd_\tau$ converges to $\cp_\tau$.}

Due to the obvious time invariance of the distributions of $\pd_\tau$ and $\cp_\tau$, it suffices to take $\tau=0$.

Let us recall from~\cite{kn:finr} that $(\ch,\dch)$ is the Hausdorff metric space generated by 
$(\Pi,d)$, and that $\yd$ converges to $\bcw$ in distribution in $(\ch,\dch)$. Resorting to 
Skorohod representation, we may take this convergence to be almost sure.
From that it follows fairly readily that 
\begin{equation}\label{conv_bd}
\dc(\pd,\cp)\to0
\end{equation}
almost surely as $\d\to0$, where $\dc$ is the Hausdorff metric associated to 
$(\cn,\bd)$, where 
\begin{equation}\label{co}
\cn:=\{((f_1,t_1),(f_2,t_2))\in\Pi'\times\Pi':\,f_2(t)\geq f_1(t)\mbox{ for }t\geq t_1\vee t_2\},
\end{equation}
the subset of ordered non crossing pairs of paths from $\Pi'\times\Pi'$.
In particular
\begin{equation}\label{conv_bd_z}
\dc(\pdz,\cpz)\to0
\end{equation}
almost surely as $\d\to0$.
It may be also readily checked that $\pd$ and $\cp$ are compact subsets of $(\cn,\bd)$.

We will strengthen~(\ref{conv_bd_z}) to
\begin{equation}\label{conv_hdh}
\hdh(\pdz,\cpz)\to0
\end{equation}
in probability as $\d\to0$.

Below we will occasionally write $\cp^{(0)}_0$ for $\cpz$.

\subsection{Argument for~(\ref{conv_hdh})}

For $\eps>0$, let $\ctd$ denote the set of pairs of neighboring points of $\d\Z\times\{0\}$ such that
the pair of paths of $\pdz$ starting from those points take longer than $\eps$ to coalesce.
Let us call such pairs of points {\em double pairs} (of points), and denote the collection of such 
paths by $\pdze$.

Let us recall that a $(0,2)$-point, or a {\em (Brownian) double point}, from $\bcw$ is a space-time point 
$(x,t)$ touched by no path of $\bcw$ starting before $t$, and from which there start (exactly) 2 paths 
from $\bcw$. See Definition 3.8 in~\cite{kn:finrp}. Let $\T$ denote the set of Brownian double points in 
$\R\times\{0\}$ such that the pairs of paths of $\cp$ starting from those points take longer than $\eps$ to 
coalesce. Let us denote the collection of such paths by $\pze$, and also occasionally by $\pdzz$.

Let $\td$, $T$ denote the cardinalities of $\ctd$, 
$\T$, 
respectively.
It follows from what we know about $\bcw$ --- see Proposition 4.3 in~\cite{kn:finr} --- that, for each 
$\eps>0$, $T$ is a finite random variable, and $\td$, $\d>0$, is a tight family of random variables.

We will call the pairs in $\pze$ and $\pdze$,
$\d\geq0$, 
as {\em (Brownian)} and {\em (discrete) double pairs} of paths (of $\cpz$ and $\pdz$, respectively), 
or {\em (Brownian)} and {\em (discrete) double paths}, for short.

\subsubsection{Convergence of double paths}\label{ss:dp}

It follows readily from~(\ref{conv_bd_z}) that there almost surely exists a 1 to 1 correspondence for each 
$\d$ small enough between discrete and Brownian double paths so that to each Brownian double path there
converges in $\bd$ distance the corresponding discrete double path as $\d\to0$. 

Now let us fix a Brownian double path, and argue that the corresponding discrete double path converges
to that Brownian path also in the $\hd$ distance. For that, it remains to show that the discrete standard 
pod associated to the discrete double path converges to that for the Brownian one.

It is enough to show that the coalescence times converge suitably. One may readily check that~(\ref{conv_bd_z})
implies that the liminf as $\d\to0$ of the discrete coalescence time (the coalescence time of the discrete 
double path) must equal or exceed the Brownian one. Let us show that the latter possibility is a null event.

If the limsup as $\d\to0$ of the discrete coalescence time exceeds the Brownian one, then given the 
constituting paths of the Brownian double path till their coalescence time and place, say $(x_0,t_0)$,
then given $M>0$, we have that for small enough $\d>0$, we may find (macroscopic) coalescing random walk 
paths starting from a (macroscopic) interval around $\d^{-1}x_0$ of length $\frac1M\d^{-1}$ at time
$\lceil\d^{-2}t_0\rceil$ which stay within distance $\frac1M\d^{-1}$ from each other without coalescing
for a time (counting from $\lceil\d^{-2}t_0\rceil$) of length $\d^{-2}T$, with $T>0$ not depending on $M$.
Consider now the locations of the intersections of those coalescing random walk paths with 
$\R\times\{\lceil\d^{-2}(t_0+T/2)\rceil\}$. We must then be able to find a pair of such points at
distance at most $\frac1M\d^{-1}$ and such that the coalescing random walk paths starting from them
stay within distance $\frac1M\d^{-1}$ from each other without coalescing for a time 
(counting from $\lceil\d^{-2}(t_0+T/2)\rceil$) of length $\d^{-2}T/2$.

Resorting again to Proposition 4.3 in~\cite{kn:finr}, we may claim that the size of the set of intersection 
locations is tight, and we will then readily conclude the argument if we show that for a single fixed pair of 
locations at distance  at most $\frac1M\d^{-1}$, the probability that the coalescing random walk paths starting 
from that pair of locations stay within distance $\frac1M\d^{-1}$ from each other without coalescing for a time 
of length $\d^{-2}T/2$ vanishes as $M\to\infty$. But this follows from standard results about the tail of the 
hitting time of the exterior of an interval by a single random walk starting from the interior of that interval.

\subsubsection{Argument for the other pairs of paths}


A similar argument as that for double paths in Subsubsection~\ref{ss:dp} shows that pairs of paths from $\pdz$ 
starting from points from different double pairs also converge to the respective Brownian pair
of paths from $\bcw$ in the $\hd$ distance.


Let us enumerate the double pairs of $\ctd=\{(x_j,y_j),\,j=1,\ldots,k\}$ in increasing order.
Let us include in $\ctd$. in case they are not already present, 
the double pairs $(x_0,y_0)$ and $(x_{k+1},y_{k+1})$ such that
$x_0\leq-1$ and $y_{k+1}\geq1$ \footnote{The paths starting from either $x_0,y_0,x_{k+1}$ or $y_{k+1}$
may not belong to $\byd$.}. Let $\l_j$ and $\rho_j$ be the portions in side $\R\times[0,\eps]$ of the
paths of $\yd$ starting from $x_j$ and $y_j$, respectively, and, for $j=0,\ldots,k$, let $\G_j$ be the open 
subset of $\R\times(0,\eps)$ bounded by 
$\rho_j$, and $\l_{j+1}$, respectively.
We will refer below to $\rho_j$ and $\l_{j+1}$, $j\geq0$, as the {\em bounding subpaths} of $\G_j$, $j\geq0$.
For $j\geq0$, let $\bar\G_j$ denote the closure of $\G_j$ and make
%
%
$\ci_j=(\R\times\{0\})\cap\bar\G_j$.

\begin{rmk}\label{rmk:coale}\mbox{}

 We note that each path starting from $\ci_j$, $j\geq1$, coalesces 
	within time $\eps$ with a path from a double path of $\pdze$, or with $\rho_0$, or wth $\l_{k+1}$ 
	--- indeed it coalesces with one of the paths bounding $\G_j$.

\end{rmk}

We next derive a bound on the {\em (horizontal) diameters} of the regions 
$\G_j^{(\d)}$, $j\geq0$, meaning the sup distance between the
bounding subpaths of $\G_j^{(\d)}$.
We will show that those diameters are bounded above
by $\eps^{1/3}$ uniformly in $j$ {\em with high probability}.

\paragraph{Bound on the diameter of $\G$ regions.}\label{bgr}

Let us start by partitioning $[-2,2]\times\{0\}$ into intervals of equal length $\eps^2$, 
and consider the event that such an interval contains more than 1 double pair. By the version 
of Proposition 4.1 in~\cite{kn:finr} which holds for $\yd$ (see argument in~\cite{kn:finra}),
together with standard facts about hitting times of the simple symmetric random walk 
on $\Z$ \footnote{Indeed, the process that comes up is the difference of two independent
simple random walks on $\Z$, but this behaves the same as the simple symmetric case as far as
this issue is concerned.},
this event has probability bounded above by $\eps^3$ uniformly for all $\d$ sufficiently
small. So the probability of the event of not finding such an interval in all of $[-2,2]\times\{0\}$ 
is bounded below by $1-\eps$ (uniformly for all $\d$ sufficiently small). We will assume 
from now on that we are in such an event.

Let us now enumerate the intervals in this partition which contain a double pair, from
left to right: $J_1,\,J_2,\ldots$. Since the right path from the double path of starting
in $J_\ell$ has to coalesce with the left path from the one in $J_{\ell+1}$, the probability
that we find a pair $J_\ell,J_{\ell+1}$ such that the respective double pairs are further than
$\eps^{1/3}$ apart is bounded above by $4\eps^{-2}$ times the probability that the
a simple random walk on $\Z$ \footnote{Idem.} starting from $\eps^{1/3}\d^{-1}$ hits the 
origin before $\eps\d^{-2}$ steps. The latter probability is bounded above by 
$e^{-\frac13\eps^{-1/3}}$ for $\eps$ small enough, uniformly in $\d$ sufficiently small. 
We thus conclude that outside an event of probability at most $e^{-\frac14\eps^{-1/3}}$, 
neighboring double pairs are at most $\eps^{1/3}$ apart 
(for $\eps$ small enough, uniformly in $\d$ sufficiently small).
A similar argument shows that outside an event of probability at most 
constant times $\eps$ for $\eps$ small enough, uniformly in $\d$ sufficiently small,
any $\G_j^{(\d)}$ region has diameter at most $\eps^{1/3}$ (here including the event 
whose probability was estimated in the previous paragraph).


\paragraph{The Brownian case.}


The claims argued above for $\pdz$ go through with no essential modification for the Brownian 
case of $\cpz$, namely, outside an event of vanishing probability as $\eps\to0$, the 
$\G$ regions have diameters at most $\eps^{1/3}$.

\subsubsection{Conclusion}

In order to conclude our argument for the validity of~(\ref{conv_hdh}), it is enough
to show that for all fixed $\eta>0$, we have that outside an event of vanishing probability as 
$\eps\to0$ uniformly in $\d$ sufficiently small 
\footnote{A $\limsup_{\d\to0}$ is taken first; and then a $\limsup_{\eps\to0}$.}, 
given a pair of paths in $\pdz$ with standard pod $p$, 
we may find a pair of paths in $\cpz$ with standard pod $q$ such that 
\begin{equation}\label{dbd}
\D(p,q)\leq\eta,
\end{equation}
and conversely, exchanging the roles of $\cpz$ and $\pdz$. The converse case is similar
to the first case, so let us argue only the first case.

It readily follows from our arguments above that we have the bound~(\ref{dbd}) for pairs of paths
from $\pdz$ starting from space-time points belonging to (not necessarily the same) double pairs; 
since the number of such double pairs is tight in $\d$ as $\d\to0$, we get~(\ref{dbd}) for all
these pairs of paths simultaneously with high probability.

In order to extend this to all paths of $\pdz$, we argue as follows. Given a pair of
paths $[(f^{(\d)}_1,0),(f^{(\d)}_2,0)]$ from $\pdz$, with standard pod $p^{(\d)}$, let 
$[(\tf^{(\d)}_1,0),(\tf^{(\d)}_2,0)]$ be the pair from $\pdz$ such that
$(\tf^{(\d)}_m,0)$ and $(f^{(\d)}_m,0)$ coalesce as described in 
Remark~\ref{rmk:coale} for $m=1,2$.

Let $q^{(\d)}$ be the standard pod of $[(\tf^{(\d)}_1,0),(\tf^{(\d)}_2,0)]$. 
From what we just argued in the previous paragraph, there exists a pair 
$[(\tf_1,0),(\tf_2,0)]$ in $\cpz$ with standard pod $q$ such that with high probability 
the bound~(\ref{dbd}) holds for the pair of pods $q^{(\d)},q$, with $\eta/2$ replacing $\eta$.

It should be quite clear from the above considerations on the diameter of the $\G$ regions that 
\begin{equation}\label{p1}
\bd([(f^{(\d)}_1,0),(f^{(\d)}_2,0)],[(\tf^{(\d)}_1,0),(\tf^{(\d)}_2,0)])
\leq 2\eps^{1/3}
\end{equation}
with high probability.

Let us now analyse $\D(p^{(\d)},q^{(\d)})$. Let $\tp^{(\d)},\tq^{(\d)}$ be the pods associted 
respectively to $p^{(\d)},q^{(\d)}$, so that the latter standard pods are obtained from the former
pods according to~(\ref{p}). Notice that $\tp^{(\d)}$ and $\tq^{(\d)}$ have the same beginning,
namely $0$. 
Let us denote their respective ends by $t^c$ and $s^c$. 

\begin{rmk}\label{ends}\mbox{}
	
	\begin{enumerate}
		\item If
		each path determining $p^{(\d)}$ is such that its portion from times $0$ to $\eps$ is 
		a bounding path of $\G_j$ for some $j$, then $p^{(\d)}=q^{(\d)}$ and 
		$\D(p^{(\d)},q^{(\d)})$ vanishes.
		\item Excluding the situation in item 1, 
		\begin{enumerate}
			\item if the starting locations of the paths determining $p^{(\d)}$ both belong to the same $\ci_j$
			for some $j$, then either $t^c=s^s$, or $t^c\leq\eps$ and $s^c=0$;
			\item excluding the situation in item 2a, if the starting locations of the paths determining $p^{(\d)}$ belong to different $\ci_j$'s, then $t^c=s^s$.
		\end{enumerate}
	    \item Whenever 	$t^c=s^s\geq2\eps$, we have that $\tp^{(\d)}(x)=\tq^{(\d)}(x)$ for $x\geq\teps:=\Psi(\eps)$, and thus
	    $p^{(\d)}(x)=q^{(\d)}(x)$ for $x\geq\teps$.
	\end{enumerate}
\end{rmk}

It follows from the above considerations that given $\eta>0$, with high probability, the following events take place simultaneously.

\begin{enumerate}
	\item If $t^c\geq2\eps$, then $\D(p^{(\d)},q^{(\d)})\leq\sum_{n\geq\teps^{-1}}2^{-n}\leq\eta/2$;
	\item If $t^c<2\eps$, then replacing $\eps$ by $2\eps$ in the analysis above (started with the bounds on diameters of the $\G$ regions on page~\pageref{bgr}), we find that $s^c\leq2\eps$ and that the diameters of both $\tp^{(\d)}$ and $\tq^{(\d)}$ are bounded above by $\eps^{1/3}$. Lemma~\ref{lem:pod} now implies the same conclusion as in the previous item, namely, $\D(p^{(\d)},q^{(\d)})\leq\eta/2$ (for all $\eps$ small enough).
\end{enumerate}

The upshot is that $\D(p^{(\d)},q^{(\d)})\leq\eta/2$ in all cases with high probability, and~(\ref{conv_hdh}) is established.


\section{Examples}
\label{ex}

\setcounter{equation}{0}

In this section, we discuss the examples mentioned in the introduction, namely persistence in the voter model and the weight distribution  at the bottom of a silo/water output distribution at a level set of a river basin (modeled by coalescing random walks, as discussed in the introduction).

\subsection{\mbox{}\hspace{-.134cm}Weight/output distribution in a silo/river basin model}
\label{wod}
As discussed in the introduction, the (properly rescaled) weight distribution at a section of the bottom of the 
silo model there described is the random measure $\mud$ on $2\d\Z\cap[-1,1]$ such that 
$\mud(\{(2\d k,0)\})$
is 
the area between the rescaled (upward, in this context --- recall the description at the introduction) random walk paths from $\byd$ issuing from $((2k-1)\d,0)$ and $((2k+1)\d,0)$.
%
The output of the 
river basin model along a section of a level line can be described by the same measure.

Let $\mu$ be the random measure on $[-1,1]$ such that given $a,b\in[-1,1]$, $a\leq b$, $\mu([a,b])$ is 
the area between the leftmost and rightmost paths from $\bcw$ issuing from $(a,0)$ and $(b,0)$, respectively.
We want to show that 
\begin{equation}\label{mud}
\mud\to\mu
\end{equation}
in distribution, with the vague topology in the space of Radon measures on $[-1,1]$. 
We argue~(\ref{mud}) next.

\subsubsection{$\mud$ and $\mu$ as images of a continuous map} 

\paragraph{A restricted space.}
Let us start by considering the following subspace of $\hh$. 
\begin{eqnarray}\nn
&\hhz=\left\{K\in\hh:\mbox{ all paths from $K$ start at time $0$;}\right.&\\\nn
&\mbox{\hspace{.5cm} for every point of $[-1,1]\times\{0\}$ there is}&\\\label{hhz}
&\left.\mbox{\hspace{.5cm} a path of $K$ starting from that point}\right\}&
\end{eqnarray}

One readily checks that $\hhz$ is a closed subset of $(\hh,\hdh)$. It then follows that $(\hhz,\hdh)$ is complete
and separable. 

In order to have the convergence~(\ref{conv_hdh}) take place in~$(\hhz,\hdh)$ as well, we must
extend the definition of $\byd$ suitably. 
A simple way to do so is by incorporating paths starting
from space-time points in $\{((2j-1)\d,(2j+1)\d)\cap[-1,1]\}\times\{0\}$, $j\in\Z$, so that their portion in $\R\times[0,\d^2]$ 
coalesces linearly at time $\d^2$ with the path from $\byd$ starting from $(2j\d,0)$;
from $((2j+1)\d,0)$, we incorporate two paths, one coalescing linearly at time $\d^2$ 
with the path from $\byd$ starting from $(2j\d,0)$;
the other coalescing linearly at time $\d^2$ 
with the path from $\byd$ starting from $(2(j+1)\d,0)$. Let us denote by $\bpd_0$ the collection of pairs of paths 
thus obtained. 
It is quite clear that $\bpd_0\in\hhz$ and that the convergence~(\ref{conv_hdh}) with $\bpd_0$ replacing $\pdz$ takes place in~$(\hhz,\hdh)$ as well --- to verify the latter claim it is enough to check that $\hdh(\bpd_0,\pdz)\to0$ as $\d\to0$,
which is a straightforward matter.
 
Let us now consider the following subset of $\hhz$.
\begin{eqnarray}\label{hhhz}
&\hhhz=\left\{K\in\hhz:\mbox{ all pairs of paths of $K$ coalesce in finite time}\right\}&
\end{eqnarray}
%
$\hhhz$ is open in  $(\hhz,\hdh)$.

Let now $\cm:\hhhz\to[0,\infty)$ be such that $\cm(K)$ is the Radon measure on $[-1,1]$ such that given 
$a,b\in[-1,1]$, $a\leq b$, $\cm(K)([a,b])$ is the area between the leftmost and rightmost paths of 
$K$~\footnote{Out of all the pairs of paths in $K$.} 
issuing from $(a,0)$ and $(b,0)$, respectively. 
The existence of the leftmost and rightmost paths of $K$ is ensured by the compactness of $K$
and the non crossing property of its paths.


Let us argue that if $K,\,K_n,n\geq1$, are in $\hhhz$ and 
\begin{equation}\label{convk}
\hdh(K_n,K)\to0 \mbox{ as } n\to\infty,
\end{equation}
then 
\begin{equation}\label{convm}
\cm(K_n)\to\cm(K) \mbox{ as } n\to\infty
\end{equation}
in the vague topology, as mentioned above (see paragraph of~(\ref{mud})).

\paragraph{Continuity of $\cm$.}

Let us first argue that 
\begin{equation}\label{c0}
\cm(K_n)([-1,1])\to\cm(K)([-1,1]). 
\end{equation}
Let $[(f^{(n)}_1,0),(f^{(n)}_2,0)]$ and $[(f_1,0),(f_2,0)]$ denote the pairs of leftmost and rightmost paths 
in $K_n$ and $K$, respectively. (\ref{convk}) implies the existence of pairs of paths 
$[(\hf^{(n)}_1,0),(\hf^{(n)}_2,0)]$ in $K$ such that
\begin{equation}\label{c1}
\hd([(f^{(n)}_1,0),(f^{(n)}_2,0)],[(\hf^{(n)}_1,0),(\hf^{(n)}_2,0)])\to0.
\end{equation}
Since $K$ is compact, there exists $[(\hf_1,0),(\hf_2,0)]$ in $K$ and a subsequence $(n')$ such that
\begin{equation}\label{c2}
\hd([(\hf^{(n')}_1,0),(\hf^{(n')}_2,0)],[(\hf_1,0),(\hf_2,0)])\to0.
\end{equation}
It is clear that $(\hf_1,0),(\hf_2,0)$ start respectively from $(-1,0),(1,0)$. If we have that either
$(\hf_1,0)\ne(f_1,0)$ or $(\hf_2,0)\ne(f_2,0)$, we get a contradiction with either the extremity of 
$(\hf^{(n')}_1,0)$ or $(\hf^{(n')}_2,0)$, respectively, for $n'$ large enough.
If follows that $(\hf_1,0)=(f_1,0)$ and $(\hf_2,0)=(f_2,0)$, and thus
\begin{equation}\label{c3}
\hd([(f^{(n)}_1,0),(f^{(n)}_2,0)],[(f_1,0),(f_2,0)])\to0,
\end{equation}
and~(\ref{c0}) follows.

To conclude, we must argue that 
\begin{equation}\label{c4}
\cm(K_n)([x,y])\to\cm(K)([x,y]),
\end{equation}
where $x,y\in[-1,1]$ are continuity points of $\cm(K)$ such that $x<y$. This means that the paths from $K$
out of $(x,0)$ and  $(y,0)$ are unique. 
Now an argument similar to the one to establish~(\ref{c0}) in the previous paragraph (however dispensing with
the contradiction part, since we have the uniqueness just alluded to), we have that 
\begin{equation}\label{c5}
\hd([(g^{(n)}_1,0),(g^{(n)}_2,0)],[(g_1,0),(g_2,0)])\to0,
\end{equation}
where $(g^{(n)}_1,0)$ and $(g^{(n)}_2,0)$ are the leftmost and rightmost paths from $K_n$ issuing from 
$(x,0)$ and $(y,0)$, respectively; $(g_1,0)$ and $(g_2,0)$ are the unique paths from $K$ issuing from 
$(x,0)$ and $(y,0)$, respectively.
(\ref{c4}) follows.

\subsubsection{Conclusion} 

The convergence~(\ref{conv_hdh}) with $\bpd_0$ replacing $\pdz$ takes place in~$(\hhz,\hdh)$,
as argued at the beginning of the previous subsubsection 
--- see paragraph right above~(\ref{hhhz}) --- 
and~(\ref{convm}) then it immediately follows that
\begin{equation}\label{bud}
\bud:=\cm(\bpd_0)\to\cm(\cpz)=\mu
\end{equation}
as $\d\to0$. To get~(\ref{mud}) it is then enough to check that for any continuous function $h:[-1,1]\to\R$,
we have that 
\begin{equation}\label{bmud}
\left| \int_{-1}^1h\,d\bud - \int_{-1}^1h\,d\mud \right| \to0
\end{equation}
as $\d\to0$, which is a straightforward matter.


\subsection{Persistence in the voter model}
\label{per}

As discussed at the introduction, the persistence probability may be expressed as follows.
For $\a\in(0,1)$, 
let $L=\{0\}\times[0,\a]$~\footnote{For convenience, we assume $\a$ to be a multiple of $\d^2$.}, and 
\begin{eqnarray}\label{hhl}
&\hhl=\left\{K\in\hh:\mbox{ all paths from $K$ start from }L\right\}.&
\end{eqnarray}
$\hhl$ is a closed subset of $(\hh,\hdh)$, and so $(\hhl,\hdh)$ is complete and separable. 
 

Let now $\cs:\hhl\to[0,\infty)$ such that for $K\in\hhl$
\begin{equation}\label{S}
\cs(K)=\sup\{t^c([(f_1,t_1),(f_2,t_2)]) :\,[(f_1,t_1),(f_2,t_2)]\in K\},
\end{equation}
that is, $\cs(K)$ is the sup of the coalescence times of pairs of paths from $K$.

$\cs$ is 
continuous in~$(\hhl,\hdh)$. 

Let us show next that
\begin{equation}\label{per1}
\cs(\pdl)\to\cs(\cpl)
\end{equation}
in distribution as $\d\to0$, where $\pdl$ and $\cpl$ are the sets of pairs of paths 
from $\pd$ and $\cp$, respectively, that start both from $L$.

\paragraph{Argument for~(\ref{per1}).}

We again assume~(\ref{conv_bd}), from which it follows that
\begin{equation}\label{conv_bd_l}
\dc(\pdl,\cpl)\to0
\end{equation}
almost surely as $\d\to0$. 

From~(\ref{conv_hdh}), we may assume that
\begin{equation}\label{conv_hdh_a}
\hdh(\pd_\a,\cp_\a)\to0
\end{equation}
almost surely as $\d\to0$. 

Let $x^{(\d)}_\ell$ and $x^{(\d)}_r$, resp.~$x_\ell$ and $x_r$, be the leftmost and rightmost points of $\R\times\{\a\}$ touched by paths from $\pdl$, resp.~$\cpl$. Theorem 3.14 from~\cite{kn:finrp} ensures that $(x_\ell,\a)$ and $(x_r,\a)$ have 
almost surely both exactly one path from $\bcw$ issuing from each of them.
Let $(f^{(\d)}_\ell,\a),(f^{(\d)}_r,\a)$ be the paths from $\byd$ issuing from $(x^{(\d)}_\ell,\a),(x^{(\d)}_r,\a)$,
respectively; and let $(f_\ell,\a),(f_r,\a)$ be the paths from $\bcw$ issuing from $(x_\ell,\a),(x_r,\a)$,
respectively

We note that 
\begin{equation}\label{coals}
\cs(\pdl)=t^c([(f^{(\d)}_\ell,\a),(f^{(\d)}_r,\a)]),\,\,\, 
\cs(\cpl)=t^c([(f_\ell,\a),(f_r,\a)]).
\end{equation}

(\ref{conv_bd_l}) now implies that 
\begin{equation}\label{xtox}
x^{(\d)}_\ell\to x_\ell,\quad x^{(\d)}_r\to x_r,
\end{equation}
and~(\ref{per1}) follows from~(\ref{conv_hdh_a}) 
\footnote{We should in this application have~(\ref{conv_hdh_a}) for paths starting from $[-r,r]\times\{0\}$ for arbitrary fixed $r$, but this of course follows as in the special case $r=1$.}
and~(\ref{xtox}) via an argument similar to one used for the 
continuity of $\cm$ in the previous subsection (see paragraph of~(\ref{c4}) above).


\section{Final comments}
\label{fin}

\setcounter{equation}{0}

A glaring issue not addressed in Section~\ref{rw} is the weak convergence of the full family of random walk paths $\pd$ to $\cp$.
The difficulty here comes in trying to obtain a uniform bound for $\D(p,q)$ over pairs of rescaled random walk paths with different pod lengths. This difficulty does not occur for pairs of paths which start at the same time and coalesce at the same time (and thus have the same pod lengths), and this facilitated the argument for~(\ref{conv_hdh}) above. The point is that in the former case we do not have the validity of something like the last point of the third item in Remark~\ref{ends}; this led in one of our attempts to the need of controlling the modulus of continuity of the {\em inverse} of a random walk path (well) after time $0$ and before hitting the origin, which we did not find a way of accomplishing.

Even if we were able to apply~(\ref{conv_hdh}) to the motivating examples, there are other interesting examples which conceivably do not follow from that result, and would require the full (or a fuller) convergence result, possibly under another, more flexible metric, in a suitable space\footnote{In a fuller, more flexible setting, one might conceivably also be able to treat the 
motivating examples 
more directly.}. One such example is the curve forming the boundary between the downward and upward families of coalescing random walk paths (those families were described in the discussion of the silo model). The convergence to the appropriately defined curve defined for the Brownian web (and its dual web) was undertaken in~\cite{kn:nr} by a direct approach. 

Another example is the full weight distribution in the silo (the collection of weights supported by all beads in the silo, not just the ones at the bottom). An intriguing point in this respect would be how to best describe this distribution (it is {\em not} a measure in the space-time plane), and in which metric space.

\bigskip

\noindent{\bf
	Acknowledgements.} At various times in the last years, besides Chuck Newman, I discussed the issues of this article with a number of colleagues and collaborators, among whom I would like to thankfully mention Rafael Grisi, who participated in early stages of this project, Pablo Ferrari, Glauco Valle, and Ant\^onio Luiz Pereira.


%

%

\end{document}